\documentclass[twocolumn,a4paper]{ifacconf}

\usepackage[latin1]{inputenc}
\usepackage{times}
\usepackage{xspace}
\usepackage[frenchb,english]{babel}
\usepackage{amsmath,amsfonts,amssymb}
\usepackage{mathrsfs}                            
\DeclareMathAlphabet\PazoBB{U}{fplmbb}{m}{n}
\renewcommand{\mathbb}[1]{\PazoBB{#1}}
\usepackage{bm}
\usepackage{enumitem}
\usepackage{xcolor}
\usepackage{graphicx}
\usepackage[authoryear,sort&compress,round]{natbib} 
\newenvironment{CustomTolerances}[2]{%
  \par 
  \pretolerance=#1\relax
  \tolerance=#2\relax
}{%
  \par
}

\newenvironment{rem}{\par\textit{Remark.}\ }{\par}
\newenvironment{nota}{\par\textit{Notations.}\ }{\par}
\newcommand{\cur}[2]{\mskip #1mu {#2} \mskip #1mu}
\newcommand{\secref}[1]{Section~\ref{#1}}
\newcommand{\thmref}[1]{Theorem~\ref{#1}}

\makeatletter
\def\@preabstractskip{8\p@ \@plus 2\p@ \@minus 2\p@} % +4\p@
\def\@belowkwskip{6\p@ \@plus 2\p@ \@minus 2\p@}
\def\@belowfmskip{28\p@} % +14\p@
\makeatother

\DeclareMathOperator*{\argmax}{argmax}
\DeclareMathOperator{\diverg}{div}

\DeclareBoldMathCommand{\EE}{E}
\DeclareBoldMathCommand{\PP}{P}

\DeclareBoldMathCommand{\fVF}{\mathrm{f}}
\DeclareBoldMathCommand{\gVF}{\mathrm{g}}
\DeclareBoldMathCommand{\jVF}{\mathrm{j}}

\newcommand    \Acal   {\mathcal{A}}
\newcommand    \Ecal   {\mathcal{E}}
\newcommand    \Filt   {\mathcal{F}}
\renewcommand  \Qset   {\mathcal{Q}}
\newcommand    \Qdisc  {{\Qset}^{\mathrm{d}}}
\newcommand    \Edisc  {{\Eset}^{\mathrm{d}}}

\newcommand    \Rplus  {\Rset_{+}}
\newcommand    \un     {\mathbb{1}}
\newcommand    \Eset   {E}
\newcommand    \EZset  {{\Eset}^0}
\newcommand    \Ecalc  {{\Ecal}_{\mathrm{c}}}
\newcommand    \Cscr   {C}
\newcommand    \Ctest  {\Cscr_{\mathrm{c}}^2{\left( \Eset \right)}}
\newcommand    \rsrc   {r^{\mathrm{src}}}
\newcommand    \Mrad   {{\mathcal{M}}_{\mathrm{c}}{\left( \Eset \right)}}

\newcommand  \ddiff {\mathrm{d}}
\newcommand  \dx    {\ddiff x}
\newcommand  \dy    {\ddiff y}
\newcommand  \dz    {\ddiff z}
\newcommand  \dH    {\ddiff H}
\newcommand  \dB    {\ddiff B}
\newcommand  \dZ    {\ddiff Z}
\newcommand  \ds    {\ddiff s}
\newcommand  \dt    {\ddiff t}

\newcommand \Gene     {L}
\newcommand \resetmap {\Psi}
\newcommand \vMes     {\mathfrak{m}}
\newcommand \sMes     {\mathfrak{s}}
\newcommand \zmin     {z_{\mathrm{min}}}
\newcommand \zmax     {z_{\mathrm{max}}}
\newcommand \dotvar   {\,\boldsymbol{\cdot}\,}
\newcommand \jin      {j^{\mathrm{in}}}
\newcommand \jout     {j^{\mathrm{out}}}

\begin{document}

\begin{frontmatter}

  % Titre (pas plus de 10 mots... raté!)
  \title{%
    A unifying formulation of the\\
    Fokker-Planck-Kolmogorov equation\\
    for general stochastic hybrid systems}

  \author{Julien Bect} 
   
  \address{Department of Signal Processing and Electronic
    Systems, Supelec, Gif-sur-Yvette, France. E-mail:
    julien.bect@supelec.fr}

  \begin{keyword}
    general stochastic hybrid systems, Markov models, continuous-time
    Markov processes, jump processes, Fokker-Planck-Kolmogorov
    equation, generalized Fokker-Planck equation
  \end{keyword}
  
  \begin{abstract}
    % Abstract of not more than 250 words.
    :\hskip 0.5em A general formulation of the Fokker-Planck-Kolmogorov
    (FPK) equation for stochastic hybrid systems is presented, within
    the framework of Generalized Stochastic Hybrid Systems (GSHS). The
    FPK equation describes the time evolution of the probability law
    of the hybrid state. Our derivation is based on the concept of
    mean jump intensity, which is related to both the usual stochastic
    intensity (in the case of spontaneous jumps) and the notion of
    probability current (in the case of forced jumps).  This work
    unifies all previously known instances of the FPK equation for
    stochastic hybrid systems, and provides GSHS practitioners with a
    tool to derive the correct evolution equation for the probability
    law of the state in any given example.
  \end{abstract}
  
\end{frontmatter}

%===== Section 1
\section{Introduction}

%===== General introduction: MP & semigroups -->GSHS
\begin{CustomTolerances}{200}{400}
  Among all continuous-time stochastic models of (nonlinear) dynamical
  systems, those with the Markov property are especially appealling
  because of their numerous nice properties. In particular, they come
  equipped with a pair of operator semigroups, the so-called backward
  and forward semigroups, which are the analytical keys to most
  practical problems involving Markov processes.  When the system is
  determined by a stochastic differential equation, these semigroups
  are generated by Partial Differential Equations (PDE) ---
  respectively the backward and forward Kolmogorov equations. The
  forward Kolmogorov PDE, also known as the Fokker-Planck equation,
  rules the time evolution~$t \mapsto \mu_t$, where $\mu_t$ is the probability
  distribution of the state~$X_t$ of the system at time~$t$. This
  paper deals with the generalization of this Fokker-Planck-Kolmogorov
  (FPK) equation to the framework of General Stochastic Hybrid Systems
  (GSHS) recently proposed by \citet{Bujorianu:2006:towards,
    Bujorianu:2004:gshs}.  
\end{CustomTolerances}

%===== Relevant literature (GSHS & generalized FPK)
The GSHS framework encompasses nearly all continuous-time Markov
models arising in practical applications, including piecewise
deterministic Markov processes \citep{Davis:1993:Markov,
  Davis:1984:piecewise} and switching diffusions
\citep{Ghosh:1997:ergodic, Ghosh:1992:switching}. Two kinds of jumps
are allowed in a GSHS: spontaneous jumps, defined by a state-dependent
stochastic intensity~$\lambda(X_t)$, and forced jumps triggered by a
so-called guard set~$G$. Generalized FPK equations have been given in
the literature, in the case of spontaneous jumps, for several classes
of models; see \citet{Gardiner:1985:handbook},
\citet{Kontorovich:1999:random-struct}, \citet{Krystul:2003:risk} and
\citet{Hespanha:2005:shs} for instance. The case of forced jumps is
harder to analyze, at the FPK level, because no stochastic intensity
exists for these jumps. Until recently, the only results available in
the literature were dealing with one-dimensional models; see
\citet{Feller1952, Feller1954} and \citet{Malhame:1985:electric}.
These results have been extended to a class of multi-dimensional
models by \citet{Bect:2006:gfk}.

%===== Contribution of the paper
The main contribution of this paper is general formulation of the FPK
equation for GSHS's. It is based on the concept of \emph{mean jump
  intensity}, which conveniently substitutes for the stochastic
intensity when the latter does not exist. This equation unifies all
previously known instances of the FPK equation for stochastic hybrid
systems, and provides GSHS practitioners with a tool to derive the
correct evolution equation for the probability law of the state in any
given example. The results presented in this paper are extracted from
the PhD thesis of the author \citep{Bect:2007:phd-thesis}.

%===== Summary
The paper is organized as follows. \secref{sec:GSHS} introduces our
notations for the GSHS formalism, together with various assumptions
that will be needed in the sequel. In \secref{sec:mean-jump-intensity}
we define the crucial concept of mean jump intensity, which is used in
\secref{sec:main-results} to derive our general formulation of the FPK
equation for GSHS's.  Section~\ref{sec:examples} concludes the paper
with a series of examples and some general remarks concerning PDEs and
integro-differential equations.

%===== Section 2
\section{General Stochastic Hybrid Systems}
\label{sec:GSHS}

\begin{CustomTolerances}{100}{300}  
  The object of interest in the GSHS formalism is a continuous-time
  strong Markov process~$X\, {=}\, (X_t)_{t \geq 0}$, with values in a
  metric space~$\EZset$. It is defined on a filtered
  space~$(\Omega,\Acal,\Filt)$, equipped with a system~$\left\{ \PP_x;\,
    x\in\EZset \right\}$ of probability measures on~$(\Omega,\Acal)$, with the
  property that~$\PP_x \left\{ X_0 = x \right\} = 1$ for all $x\in\EZset$
  (i.e., $X$ starts from~$x$ under~$\PP_x$). As usual, $\EE_x$ denotes
  the expectation operator corresponding to~$\PP_x$.
\end{CustomTolerances}

It is assumed that, for each~$\omega \in \Omega$, the samplepath $t \mapsto
X_t(\omega)$ is right-continuous, has left limits~$X_t^-(\omega)$ in the
completion~$\Eset$ of~$\EZset$, and has a finite number of jumps,
denoted by~$N_t(\omega)$, on the interval~$[0;t]$ for all~$t \geq 0$. The
last condition can be seen as a ``pathwise non-Zenoness'' requirement.
We will denote by~$\tau_k$ the $k^{\text{th}}$ jump time, with $\tau_k
\,{=}\, {+}\, \infty$ if there is less than $k$~jumps.

%===== Subsection 2.1
\subsection{The hybrid state space}
\label{subsec:hybrid-ss}

The (completed) state-space of the model is assumed to have a hybrid
structure: $\Eset = \cup_{q\in\Qset}\, \{q\} × \Eset_q$, where~$\Qset$ is a
finite or countable set, and each~$\Eset_q$ is either the closure of
some connected open subset $D_q \subset \Rset^{n_q}$ ($n_q \geq 1$) or any
singleton space (in which case we set $n_q = 0$). The state at
time~$t$ can therefore be written as a pair $X_t \,{=}\, (Q_t, Z_t)$,
where~$Q_t \in \Qset$ and $Z_t \in \Eset_{Q_t}$. We denote by $\Qdisc =
\left\{ q \in \Qset \bigm| n_q = 0 \right\}$ the set of all ``purely
discrete'' modes, and by $\Edisc = \cup_{q\in\Qdisc}\, \{q\} × \Eset_q$ the
corresponding subset of~$\Eset$.

The state space~$\Eset$ is regarded as the disjoint sum of the
sets~$\Eset_q$, $q \in \Qset$, and endowed with the disjoint union
topology\footnote{which is (here) locally compact, separable and
  completely metrizable}. We denote by~$\Ecal$ the Borel $\sigma$-algebra,
and by~$\Ecalc$ the subsets of all relatively compact~$\Gamma \cur{4}{\in}
\Ecal$.  Moreover, we define a ``volume measure'' on~$\Eset$ by the
relation
\begin{equation*}
  \vMes(\Gamma) = \sum_{q\not\in\Qdisc} \vMes_q(\Gamma \cap \Eset_q)
  + \sum_{x \in \Edisc} \delta_x(\Gamma) \,, \quad \Gamma \in \Ecal,
\end{equation*}
where~$\vMes_q$ is the $n_q$-dimensional Lebesgue measure on~$\Eset_q$
and~$\delta_x$ the Dirac mass at~$x$. (Note that~$\Eset_q \subset
\Rset^{n_q}$ has been tacitly identified with $\{q\} × \Eset_q
\subset \Eset$.)

Let $\partial\Eset_q$ be the boundary of~$\Eset_q$ in~$\Rset^{n_q}$, with the
convention that~$\partial\Eset_q = \varnothing$ when~$n_q = 0$. We define the
boundary~$\partial\Eset$ of the state space by the relation $\partial\Eset =
\cup_{q\in\Qset}\, \{q\} × \partial\Eset_q$, and the \emph{guard set} by $G = \Eset
\setminus \EZset$. It is \emph{not} required that $G = \partial\Eset$.

% We assume that:
% \begin{assum} 
%   \label{ass:ss}
%   \begin{enumerate}
%   \item \label{ass:ss1} $G = \partial\Eset$,
%   \item \label{ass:ss1} $\partial\Eset_q$ is of class~$\Cscr^2$, for each~$q$
%     such that $n_q \geq 2$.
%   \end{enumerate}
% \end{assum}

% \begin{rem}
%   Both these assumptions are specific to this paper, not to the GSHS
%   formalism. They can be relaxed in some cases: for instance, the
%   boundary could be allowed to have corners, as long as some version
%   of the divergence theorem remains available.
% \end{rem}

\begin{nota}
  Let $\mu:\Ecal \to \Rset$ be a (signed) measure, $K: \Eset
  × \Ecal \mapsto \Rset$ a kernel and $\varphi: \Eset \to
  \Rset$ a measurable function. The following notations will be used
  throughout the paper, assuming the integrals exist: $(\mu K)(\dy) =
  \int \mu(\dx)\, K(x,\dy)$, $(K \varphi)(x) = \int K(x,\dy)\, \varphi(y)$ and $\mu
  \varphi = \int \mu(\dx)\, \varphi(x)$.
%   \begin{gather*}
%     \mu \varphi = \int \mu(\dx)\, \varphi(x) \,, 
%     \quad (\mu K)(\dy) = \int \mu(\dx)\, K(x,\dy)
%     \,, \\
%     \text{and} \quad (K \varphi)(x) = \int K(x,\dy)\, \varphi(y) \,.
%   \end{gather*}
\end{nota}

%===== Subsection 2.2
\subsection{A stochastic differential equation with jumps}
\label{subsec:SDE-with-jumps}

A vector field~$\gVF$ on~$\Eset$ is regarded as a first order
differential operator with respect to the continuous variables: its
action on a continuously differentiable function $\varphi \in \Cscr^1(\Eset)$
will be denoted by $\gVF \varphi$, where $(\gVF \varphi)(q,z) \cur{5}{=}\!
\sum_{i=1}^{n_q} \gVF^i(q,z) \frac{\ddiff\varphi}{\dz^i}(q,z)$ on $\Eset \setminus
\Edisc$ and $\gVF \varphi \cur{5}{=} 0$ on~$\Edisc$. The number of
``components'' of~$\gVF$ depends on the mode~$q$: to simplify the
notations, we shall agree that the indexes~$i$ and~$j$ always
correspond to summations on the number of continuous variables, and
drop the explicit dependence on~$q$. For instance, the definition
of~$\gVF \varphi$ can be rewritten as $\gVF \varphi = \sum_i \gVF^i \frac{\partial\varphi}{\partial
  z^i}$.

The process~$X$ is assumed to be driven by an Itô stochastic
differential equation between its jumps: there exist $r+1$ smooth
vector fields $\fVF_l$ and a $r$-dimensional Wiener process~$B$ such
that, in mode~$q\in\Qset \setminus \Qdisc$,
\begin{equation}
  \label{eq:SDE-Z}
  \dZ_t \;=\; \fVF_0(q,Z_t)\, \dt
  + \sum_{l=1}^r \fVF_l (q,Z_t)\, \dB^l_t \,.
\end{equation}
In other words, for all $\varphi \in \Cscr^2(\Eset)$, $X$ satisfies
the following generalized Itô formula
\begin{equation*}
  \label{eq:SDE-X}
  \begin{split}
    \varphi(X_t) - \varphi(X_0) =\; & \int_0^t (\Gene \varphi)(X_s) \,\ds
    + \sum_{l=1}^r \int_0^t (\fVF_l \varphi)(X_s) \,\dB^l_s \\
    & + \sum_{0 < \tau_k \leq t} \bigl( 
    \varphi(X_{\tau_k}) - \varphi(X_{\tau_k}^-) \bigr)\,,
  \end{split}
\end{equation*}
where $\Gene$ is the differential generator associated
with~(\ref{eq:SDE-Z}), i.e. $\Gene = \sum_i \fVF_0^i\, \frac{\partial}{\partial z^i} + \tfrac{1}{2}
\sum_{i,j}\, \left( \sum_{l=1}^r \fVF_l^i \fVF_l^j \right) \frac{\partial^2}{\partial z^i
  \partial z^j}$. We make the following smoothness assumptions:
\begin{assum}
  The drift~$\fVF_0$ is of class~$\Cscr^1$, and the other vector
  fields~$\fVF_l$, $1 \leq l \leq r$, are of class~$\Cscr^2$.
\end{assum}

%===== Subsection 2.3
\subsection{Two different kinds of jumps}

We assume that there exists a Markov kernel~$K$ from~$E$ to $\EZset$
and a measurable locally bounded function $\lambda:\EZset \to
\Rplus$, such that the following \emph{Lévy system identity} holds for
all $x \in \EZset$, $t \geq 0$, and for all measurable $\varphi:\Eset ×
\EZset \to \Rplus$:
\begin{equation*}
  \EE_x \left\{ \vphantom{\int_0^t} \sum\nolimits_{0 < \tau_k \leq t} \varphi(X_{\tau_k}^-,X_{\tau_k})
  \right\} = \EE_x \left\{ \int_0^t (K\varphi)(X_s^-)\, \dH_s \right\}
\end{equation*}
where $(K\varphi)(y) = \int_{\EZset} K(y,\dy')\, \varphi(y,y')$ and $H$ is the
predictable increasing process defined by
\begin{equation}
  \label{eq:2}
  H_t = \int_0^t \lambda(X_s)\, \ds + \sum_{\tau_k\leq t} \un_{X_{\tau_k}^- \in G} \,.
\end{equation}
The first part corresponds to \emph{spontaneous} jumps, triggered
``randomly in time'' with a stochastic intensity $\lambda(X_t)$, while the
other part corresponds to \emph{forced} jumps, triggered when~$X$ hits
the guard set~$G$.

\begin{rem}
  The terms ``spontaneous'' and ``forced'' seem to have been coined by
  \citet{Bujorianu:2003:modelling-framework}. They are closely related
  to the probabilistic notions of predictability and total
  inaccessibility for stopping times \citep[see, e.g.,][chapter~VI,
  §§12--18]{RogersWilliams:vol-2}, but be shall not discuss this point
  further in this paper.
\end{rem}

\begin{rem}
  The pair $(K,H)$ is a \emph{Lévy system} for the process~$X$ in the
  sense of \citet[definition~6.1]{Walsh:1972:terminal-times}. Most
  authors require that~$H$ be continuous in the definition of a Lévy
  system, thereby disallowing predictable jumps.
\end{rem}

%====== Section 3
\section{Mean jump intensity}
\label{sec:mean-jump-intensity}

From now on, we assume that some initial probability law~$\mu_0$ has
been chosen, with $\mu_0(G)=0$ since the process cannot start from~$G$.
All expectations will be taken, without further mention, with respect
to the probability $\PP_{\mu_0} = \int \mu_0(\dx) \PP_x$.

%====== Subsection 3.1
\subsection{Definition and link with the usual stochastic intensity}
\label{subsec:def-and-link}

It is assumed from now on that $\EE(N_t) \,{<}\, {+}\,\infty$. This is a
usual requirement for stochastic hybrid processes\footnote{See, e.g.,
  \citet{Davis:1984:piecewise} or \citet{Bujorianu:2004:gshs}.}, which
is clearly stronger than piecewise-continuity of the samplepaths. Its
being satisfied depends not only on the dynamics of the system but
also on the initial probability law~$\mu_0$.

In order to introduce the main concept of this section, let us define
a (positive, unbounded) measure~$R$ on $\Eset × \left( 0;+\infty \right)$
by
\begin{equation*}
  R\left( \mskip -2mu A \right) 
  = \EE_{\mu_0} \biggl\{ \sum\nolimits_{k\geq1}\, \un_{A} \left( X_{\tau_k}^-, \tau_k 
  \right) \biggr\} \,.
\end{equation*}
For any~$\Gamma \,{\in}\, \Ecal$, the quantity $R\left( \Gamma × (0;t] \right)$ is
the expected number of jumps starting from~$\Gamma$ during the time
interval~$(0;t]$.
\begin{defn}
  Suppose that there exists a mapping~$r:t \mapsto r_t$, from~$[0;+\infty)$
  to the set of all positive bounded measures on~$\Eset$, such that,
  for all~$\Gamma \in \Ecal$,
  \begin{enumerate}
  \item $t \mapsto r_t(\Gamma)$ is measurable,
  \item for all~$t \geq 0$, $R\left( \Gamma × (0;t] \right) = \int_0^t\,
    r_s(\Gamma)\, \ds$.
  \end{enumerate}
  Then~$r$ is called the \emph{mean jump intensity} of the process~$X$
  (started with the initial law~$\mu_0$).
\end{defn}

Let us split~$R$ into the sum of two measures~$R^0$ and~$R^G$,
corresponding respectively to the spontaneous and forced jumps of
the process. Then, using the Lévy system identity, it is easy to see
that a mean jump intensity~$r^0$ always exist for the spontaneous
part~$R^0$: it is given by
\begin{equation*}
  r^0_t(\Gamma) = \EE \bigl( \lambda(X_t)\, \un_{X_t \in \Gamma} 
  \bigr) = \int_{\Gamma} \lambda(x)\, \mu_t(\dx)\,.
\end{equation*}
In other words: for spontaneous jumps, a mean jump intensity always
exists, and it is the expectation of the stochastic jump
intensity~$\lambda(X_t)$ on the event~$\{ X_t \in \Gamma \}$.

Forced jumps are more problematic. The Lévy system identity is
powerless here, since no stochastic intensity exists (because forced
jumps are predictable). All hope is not lost, though: a simple example
will be presented in the next subsection, proving that a mean jump
intensity can exist anyway. This is fortunate, since the existence of
a mean jump intensity will be an essential ingredient for our unified
formulation of the generalized FPK equation. See
subsection~\ref{subsec:forced-ex} for further details on that issue.

%====== Subsection 3.2
\subsection{Where $\mu_0$ comes into play: an illustrative example}

Consider the following hybrid dynamics on $\Eset \,{=}\, [0;1]$: the
state~$X_t$ moves to the right at constant speed $v > 0$ as long as it
is in~$\EZset \,{=}\, [0;1)$, and jumps instantaneously to~$0$ as soon
as it hits the guard~$G=\{1\}$ (i.e., the reset kernel is such
that~$K(1,\,\cdot\,) \,{=}\, \delta_0$).

If we take $\mu_0 = \delta_0$ for the initial law, then the process jumps
from~$1$ to~$0$ each time~$t$ is a multiple of~$1/v$, i.e. $\tau_k = k/v$
and $X_{\tau_k}^- = 1$ almost surely. There is therefore no mean jump
intensity in this case, since $R = \sum_{k\geq1}\, \delta_{(1,\,k/v)}$.

Now take~$\mu_0$ to be the uniform probability on~$[0;1]$ (which is,
incidentally, the only stationary probability law of the process).
Then
\begin{align*}
  R\bigl( \Gamma × (0;t] \bigr) 
  & \;=\; \delta_1(\Gamma)\, \int_0^1\, \argmax_{k\geq1} \left\{ 
    \frac{k-x}{v} \leq t \right\}\, \dx \\
  & \;=\; \delta_1(\Gamma)\, \int_0^1\, \left\lceil vt+x \right\rceil\, \dx \\
  & \;=\; vt\, \delta_1(\Gamma) \,,
\end{align*}
where $\left\lceil vt+x \right\rceil$ is the smaller integer greater or equal to
$vt + x$. Therefore the mean jump intensity exists in this case, and
is equal to $v\, \delta_1$ (it is of course time-independent, since~$\mu_0$
is stationary). In particular, the global mean jump intensity is
$r_t(\Eset) = v$.

%==============================================================================

\section{Generalized FPK equation}
\label{sec:main-results}

%===== 4.1
\subsection{A weak form of the FPK equation}
\label{subsec:weak-FPK}

Taking expectations in \ref{eq:SDE-X}, the following \emph{generalized
  Dynkin formula} is obtained: for all compactly supported $\varphi
\in \Cscr^2(\Eset)$ and all~$t \geq0$,
\begin{equation}
  \label{eq:Dynki-with-expectations}
  \begin{split}
    \EE&\left\{ \varphi(X_t) - \varphi(X_0) \right\}
    \;=\;
    \EE \left\{ \int_0^t (L\varphi)(X_s)\, \ds \right\}\\
    &\qquad\qquad\qquad+ 
    \EE \Biggl\{  \sum_{0 < \tau_k \leq t} 
    \varphi(X_{\tau_k}) - \varphi(X_{\tau_k}^-) \Biggr\}
    \,.    
  \end{split}
\end{equation}
Let us assume the existence of a mean jump intensity~$r_t$ at all
times. Then~\eqref{eq:Dynki-with-expectations} can be rewritten as
\begin{equation}
  \label{eq:Dynkin-rewritten}
  \left( \mu_t - \mu_0 \right) \varphi \;=\;
  \int_0^t \mu_s (L\varphi)\, \ds + 
  \int_0^t r_s (K - I) \varphi\, \ds \,,
\end{equation}
where~$\mu_t$ is the law of~$X_t$ and~$I$ is the ``identity kernel''
on~$\Eset$, i.e. the kernel defined by $I(y,\dy') = \delta_y(\dy')$.
Formally differentiating~\eqref{eq:Dynkin-rewritten} yields
\begin{equation}
  \label{eq:weak-FPK}
  \mu'_t \;=\; \Gene^* \mu_t + r_t (K - I) \,,
\end{equation}
where $t \mapsto \mu'_t$ is the ``derivative'' of~$t \mapsto \mu_t$ (in a
sense to be specified later), and~$\Gene^*$ the adjoint of~$\Gene$ in
the sense of distribution theory. 

Equation~\eqref{eq:weak-FPK} begins like the usual Fokker-Planck
equation for diffusion processes ($\mu'_t = \Gene^* \mu_t$) and ends with
an additional term that accounts for the jumps of the process.

\begin{defn}
  \label{def:GFPK}
  We will say that~$t \mapsto \mu_t$ is a solution in the weak sense of the
  \emph{generalized FPK equation} for the GSHS if
  \begin{enumerate}[label=\alph*),ref=\thedefn.\alph*]
  \item \label{def:GFPK:a} there exists a mean jump intensity~$t \mapsto
    r_t$,
  \item \label{def:GFPK:b} there exists a mapping~$t \mapsto \mu'_t$,
    from~$[0;+\infty)$ to the space~$\Mrad$ of all Radon measures
    on~$\Eset$, such that $t \mapsto \mu_t(\Gamma)$ is absolutely continuous with
    a.e.-derivative~$t \mapsto \mu'_t(\Gamma)$, for all~$\Gamma \in \Ecalc$,
  \item \label{def:GFPK:c} $\Gene^* \mu_t$ is a Radon measure for all
    $t\geq 0$,
  \item \label{def:GFPK:d} equation~\eqref{eq:weak-FPK} holds as an
    equality between Radon measures, i.e. $\mu'_t(\Gamma) = (\Gene^* \mu_t)(\Gamma)
    + r_t (K-I)(\Gamma)$ for all $t \geq0$ and all~$\Gamma \in \Ecalc$.
  \end{enumerate}
\end{defn}

Such a weak form of the FPK equation is the price to pay for a unified
treatment of both kind of jumps.  Conditions~\ref{def:GFPK:a}
and~\ref{def:GFPK:b} can be seen as smoothness requirements with
respect to the time variable, and~\ref{def:GFPK:c} with respect to the
space variables.

%====== 4.2
\subsection{``Physical'' interpretation}
\label{subsec:phy-int}

The usual FPK equation admits a well-known physical interpretation as
a conservation equation for the ``probability mass'' \citep[see
e.g.][]{Gardiner:1985:handbook}. Indeed, assuming the existence of a
smooth pdf~$p \in \Cscr^{2,1}(\Eset × \Rplus)$, the
equation~$\mu'_t = \Gene^* \mu_t$ can be rewritten as a conservation
equation $\partial p_t / \partial t = \diverg ( \jVF_t )$, with the
\emph{probability current}~$\jVF_t$ defined by
\begin{equation}
  \label{eq:proba-curr}
  \jVF^i_t = \fVF_0^i\, p_t
  - \frac{1}{2} \sum_j \frac{\partial (a^{ij} p_t)}{\partial z^j}\,,
  \quad
  a^{ij} = \sum_{l=1}^r \fVF_l^i \fVF_l^j \,.
\end{equation}

The additional ``jump term'' in the generalized FPK equation, admit a
nice physical interpretation as well. To see this, let us rewrite it
as the difference of two bounded positive measure: $r_t (K-I) =
\rsrc_t - r_t$, where $\rsrc_t = r_t K$. Therefore~$r_t$ and~$\rsrc_t$
behave respectively as a \emph{sink} and a \emph{source} in the
generalized FPK equation: for each~$\Gamma\in\Ecal$, $r_t(\Gamma)\,\dt$ is
the probability mass leaving the set~$\Gamma$ during~$\dt$, because of the
jumps of the process, while~$\rsrc_t(\Gamma)\, \dt$ is the probability mass
entering~$\Gamma$. 

These two measures are in fact connected by the reset
kernel~$K(x,\dy)$. In particular, the relation~$r_t(\Eset) =
\rsrc_t(\Eset)$ holds at all times~$t \geq 0$, ensuring that the total
probability mass is conserved. Moreover, introducing the
measures~$W_t(\dx,\dy) = r_t(\dx) K(x,\dy)$, we have $r_t = \int
W(\cdot,\dx)$, $\rsrc_t = \int W(\dx,\cdot)$ and the generalized FPK
equation can be rewritten more symmetrically as
\begin{equation*}
  \mu'_t \;=\; \Gene^* \mu_t + 
  \int \left( W_t(\dx,\cdot) - W_t(\cdot,\dx) \right) \,.
\end{equation*}
It appears clearly, under this form, as a generalization of the
\emph{differential Chapman-Kolmogorov formula} of
\citet[equation~3.4.22]{Gardiner:1985:handbook} --- which only allows
spontaneous jumps.

%===== 4.3
%\subsection{Main theorem}
\subsection{Sufficient conditions for the existence of a weak solution}

The main result of this paper show that the various requirements of
definition~\ref{def:GFPK} are not independent. We denote by~$\left| \nu
\right|$ the total variation measure of a Radon measure~$\nu$, which is
finite on~$\Ecalc$. We shall say that a function $t \mapsto \nu_t$
from~$[0;\infty)$ to~$\Mrad$ is right-continuous (resp. locally integrable)
is $t \mapsto \nu_t \varphi$ is right-continuous (resp. locally integrable)
for all bounded measurable $\varphi:\Eset \to \Rset$.

\begin{thm}
  \label{THM:MAIN}
  Consider the following assumptions:
  \begin{enumerate}[label=\alph*),ref=\thedefn.\alph*]
  \item\label{THM:MAIN:a} there exists a mean jump intensity~$r$
    (\ref{def:GFPK:a}), such that~$t \mapsto r_t$ is right-continuous,
  \item\label{THM:MAIN:b} $t \mapsto \mu_t$ is differentiable in the sense
    of~\ref{def:GFPK:b}, $t \mapsto \mu'_t$ is right-continuous and $t \mapsto
    \left| \mu'_t \right|$ locally integrable,
  \item\label{THM:MAIN:c} $\Gene^* \mu_t$ is a Radon measure for all $t\geq
    0$ (\ref{def:GFPK:c}), $t \mapsto \Gene^* \mu_t$ is right-continuous and
    $t \mapsto \left| \Gene^* \mu_t \right|$ is locally integrable.
  \end{enumerate}
  If any two of these assumptions hold, then the third holds as well
  and $t \mapsto \mu_t$ is a solution in the weak sense of the
  generalized FPK equation.
\end{thm}

\newcommand \Ha {\ref{THM:MAIN:a}\xspace}
\newcommand \Hb {\ref{THM:MAIN:b}\xspace}
\newcommand \Hc {\ref{THM:MAIN:c}\xspace}

\begin{CustomTolerances}{700}{2000}
  The proof of this theorem is given in appendix~\ref{sec:proof}. We
  will not try to give general conditions under which
  assumptions~\ref{THM:MAIN:a}--\ref{THM:MAIN:c} are satisfied, since
  such conditions would inevitably be, in the general setting of this
  paper, very complicated (involving the initial law~$\mu_0$, the vector
  fields~$\gVF$ of the stochastic differential equation, the geometry
  of the state space~$\Eset$ and the reset kernel~$K$).
\end{CustomTolerances}

%====== 4.4
\subsection{The case when a piecewise smooth pdf exists}
\label{subsec:piecewise-smooth}

Equation~\eqref{eq:weak-FPK} is an evolution equation for the
measure-valued function~$t \mapsto \mu_t$. In most situations of
practical interest, the measures~$\mu_t$ admit a pdf~$p_t$, with respect
to the volume measure~$\vMes$ on~$\Eset$ (sometimes with an additional
singular measure, like a linear combination of Dirac masses, but this
case will not be discussed here). If the function~$p:(x,t)\mapsto
p_t(x)$ is smooth enough, at least piecewise, then
equation~\eqref{eq:weak-FPK} simultaneously gives birth to an
evolution equation for~$t \mapsto p_t$ and to static relations that
hold for all~$t \geq 0$ (so-called ``boundary conditions'', although
the name is not entirely appropriate here). This can be done quite
generally, using some additional measure-theoretic tools for which
there is no room in this paper. The reader is referred to
\citet[§IV.2.C]{Bect:2007:phd-thesis} for more on this issue.

%======= 5
\section{Examples}
\label{sec:examples}

%======= 5.1
\subsection{A class of models with spontaneous jumps}
\label{subsec:spont}

Our first series of examples covers a large family of models without
forced jumps ($G = \varnothing$). The reset kernel~$K$ is assumed to
satisfy the following assumption:
\begin{assum}
  \label{assumpt:dual-kernel}
  There exists a kernel~$K^*$ on~$\Eset$ such that
  \begin{equation*}
    \vMes(\dx)\, K(x,\dy) = \vMes(\dy)\, K^*(y,\dx) \,.
  \end{equation*}
\end{assum}
\removelastskip
(We do \emph{not} assume that~$K^*$ is a Markov kernel, i.e.  that
$K^*(y,\cdot)$ is a probability measure for all~$y$.) The following result
is an easy consequence of \thmref{THM:MAIN}:
\begin{cor}
  \label{cor:fpk-spont}
  If there exists a pdf~$p \in \Cscr^{2,1}(\Eset × \Rplus)$, then
  the measures~$r_t$ and~$\rsrc_t$ are absolutely continuous with
  respect to~$\vMes$,
  \begin{equation*}
    \frac{\ddiff r_t}{\ddiff \vMes} \,=\, \lambda\, p_t
    \,, \quad
    \frac{\ddiff \rsrc_t}{\ddiff \vMes} \,=\, K^*\left(\lambda\, p_t\right)
    \,,
  \end{equation*}
  and the following evolution equation holds:
  \begin{equation}
    \label{eq:FPK-spont}
    \frac{\partial p_t}{\partial t}
    \;=\;
    \Gene^* p_t \,+\, 
    K^*\left(\lambda\, p_t\right) \,-\, \lambda\, p_t \,.
  \end{equation}
\end{cor}
Assumption~\ref{assumpt:dual-kernel} holds for several classes of
models known in the literature: pure jump processes with an absolutely
continuous reset kernel, the switching diffusions of
\citet{Ghosh:1997:ergodic, Ghosh:1992:switching} and also the SHS of
\citet{Hespanha:2005:shs}.

\begin{exmp} \label{ex:pure-jumps} Pure jump processes occur when
  $\Gene = 0$, i.e. when there is no continuous dynamics. We consider
  here the case where~$K$ is absolutely continuous: $K(x,\dy) =
  k(x,y)\, \vMes(\dy)$. For instance, if the amplitude of the jumps is
  independent of the pre-jump state and distributed the pdf~$\rho$, then
  $k(x,y) = \rho(y-x)$.  In this case
  Assumption~\ref{assumpt:dual-kernel} holds with $K^*(x,\dy) =
  k(y,x)\, \vMes(\dy)$. Introducing the function $\gamma(x,y) = \lambda(x)
  k(x,y)$, equation~\ref{eq:FPK-spont} turns into the well-known
  \emph{master equation} \citep[eq.~3.5.2]{Gardiner:1985:handbook}:
  \begin{equation*}
    \frac{\partial p}{\partial t}(y,t)
    \;=\; \int \bigl( 
    \gamma(x,y) p(x,t) - \gamma(y,x) p(y,t) \bigr)\, 
    \vMes(\dx) \,.
  \end{equation*}
  In particular, when all modes are purely discrete ($n_q = 0$), this
  is just the usual forward Kolmogorov equation for a continuous-time
  Markov chain.
\end{exmp}

\begin{exmp} \label{ex:switching-diff} In the case of switching
  diffusions, the state space is of the form~$\Eset = \Qset × \Rset^n$
  (with~$\Qset$ a countable set and~$n\geq 1$) and the reset kernel of
  the form
  \begin{equation*}
    K\bigl((q,z),\cdot\bigr) = \sum_{q'\neq q}\, \pi_{qq'}(z)\, \delta_{(q',z)} \,,
  \end{equation*}
  where $\pi(z) = \left( \pi_{q q'}(z) \right)$ is a stochastic matrix for
  all~$z \in \Rset^n$. Assumption~\ref{assumpt:dual-kernel} is
  fulfilled with~$K^*$ defined by
  \begin{equation*}
    K^* \bigl((q,z),\cdot\bigr) 
    = \sum_{q'\neq q}\, \pi_{q'q}(z)\, \delta_{(q',z)}
    \,.
  \end{equation*}
  Equation~\ref{cor:fpk-spont} becomes in this case the familiar
  generalized FPK equation for switching diffusion processes
  \citep[see, e.g.,][]{Kontorovich:1999:random-struct,
    Krystul:2003:risk}: for all $x=(q,z)\in \Eset$ and $t\geq0$,
  \begin{equation*}
    \frac{\partial p}{\partial t}(x,t) = (\Gene^* p_t)(x) +
    \sum_{q'\neq q} \lambda_{q'\!q}(z)\, p_t(q',z) - \lambda(x)\, p_t(x) \,,
  \end{equation*}
  where $\lambda_{q'\!q}(z) = \lambda(q',z)\, \pi_{q'\!q}(z)$.
\end{exmp}

\begin{exmp} The SHS of \citet{Hespanha:2005:shs} are also defined
  on~$\Eset = \Qset × \Rset^n$, but this time the post-jump
  state~$X_{\tau_k}$ is determined by applying a reset map~$\Psi:\Eset
  \to \EZset$ to the pre-jump state~$X^-_{\tau_k}$, $\Psi$ being
  chosen randomly in a finite of reset maps~$\Psi_k$. The reset kernel
  can therefore be written as
  \begin{equation*}
    K(x,\cdot) = \sum_k\, \pi_k(x)\, \delta_{\resetmap_k(x)} \,,
  \end{equation*}
  with~$\pi_k(x)$ the probability of choosing the reset map~$\Psi_k$ given
  that $X^-_{\tau_k} = x$. Provided that the functions $\Psi_k$ are local
  $\Cscr^1$-diffeomorphisms, the kernel~$K$
  fulfills Assumption~\ref{assumpt:dual-kernel} with
  \begin{equation*}
    K^*(x,\cdot) = \sum_k\, \sum_{y\in\resetmap_k^{-1}(\{x\})}
    \pi_k(y)\, \bigl| J_k(y) \bigr|^{-1}\, \delta_y \,,
  \end{equation*}
  where~$J_k(y)$ is the Jacobian determinant of~$\Psi_k$ at~$y$.
  Therefore, introducing a stochastic intensity~$\lambda_k = \lambda\, \varrho_k$ for
  each one of the reset maps, we recover thanks to
  Corollary~\ref{cor:fpk-spont} the generalized FPK equation given by
  \citet[p.~1364]{Hespanha:2005:shs}: 
  \begin{equation*}
    \begin{split}
      \frac{\partial p}{\partial t}(x,t) 
      & \;=\;
      (\Gene^* p_t)(x) \\
      & + \sum_k\,
      \sum_{y\in\resetmap_k^{-1} (\{x\}) }
      \left(
        \frac{ \lambda_k\, p_t}{ \left| J_k \right|}(y)\, 
        - (\lambda_k\, p_t)(x)
      \right) \,.      
    \end{split}
   \end{equation*}
\end{exmp}

%====== 5.2
\subsection{A class of models with forced jumps}
\label{subsec:forced-ex}

The measure-valued formulation of the generalized FPK
equation~equation~\eqref{eq:weak-FPK} paves the way for an easier
proof of some recent results \citep{Bect:2006:gfk}, concerning GSHS
with forced jumps and deterministic resets. A typical example of this
class of process is the thermostat model
of~\citet{Malhame:1985:electric}. Since a complete statement and proof
of these results would be too long for this paper, we shall only
provide an illustrative example. The interested reader is referred to
the PhD thesis of the author \citep[IV.2.C and
IV.3.C]{Bect:2007:phd-thesis}. A thorough treatment will appear in a
forthcoming publication.

\begin{exmp} \label{ex:thermo} Let us consider a GSHS without
  spontaneous jumps ($\lambda = 0$), whose hybrid state space is defined by
  $\Qset = \{0,1\}$, $\Eset_0 = [\zmin;+\infty)×\Rset^{n-1}$, and $\Eset_0 =
  (-\infty;\zmax]×\Rset^{n-1}$ (where $\zmin \,{<}\, \zmax$). Assume that
  the guard~$G$ is the whole boundary~$\partial\Eset$, and that the reset map
  is defined by $\Psi(q,z) = (1-q,z)$. In other words, the discrete
  component~$Q_t$ switches from~$0$ to~$1$ when $Z^1_t$ reaches the
  lower threshold~$\zmin$, and switches back to~$0$ when $Z^1_t$
  reaches the upper threshold~$\zmax$.

  For such a hybrid structure, it is easily shown using
  \thmref{THM:MAIN} that no~$\Cscr^{2,1}$ solution can exist.
  Consider the set $G' = \Psi(G)$, which is the disjoint unions of two
  ``hyperplanes'' in $\EZset$. A careful examination
  of~\eqref{eq:weak-FPK} suggests to look for solution that are of
  class~$\Cscr^{2,1}$ on~$\EZset \setminus G'$, possibly with a discontinuity
  on~$G'$. If the process effectively has a pdf~$p$ satisfying these
  assumptions, then it can be proved using \thmref{THM:MAIN} that:
  \begin{enumerate}
  \item The usual Fokker-Planck equation, $\partial p_t/\partial t =
    \Gene^* p_t$, holds on the four components of~$\EZset \setminus
    G'$,
  \item The jumps are accounted for by the static relation $\jout_t =
    \jin_t \circ \psi$ on~$G$, at all times~$t \geq0$, where~$\jout_t$
    and~$\jin_t$ are the outgoing and ingoing probability current,
    respectively defined on~$G$ and~$G'$
    (see~\eqref{eq:proba-curr} for the defintion of the probability
    current).
  \item The mean jump intensity~$r_t$ is supported by~$G$ and given by
    the outgoing flux of the probabily current~$\jVF_t$, i.e. $r_t(\Gamma)
    = \int_{\Gamma \cap G} \jout_t \ddiff\sMes$, where~$\sMes$ is the
    surface measure. 
  \item Finally, for each~$x \in G$ such that at least one of the
    ``noise driven'' vector fields~$\gVF_l$ ($1 \leq l \leq r$) is
    transverse to~$G$, the pdf has to satisfy the so-called
    \emph{absorbing boundary condition} $p_t(x) = 0$. For similar
    reasons, $p_t$ has to be continuous at each~$x \in \Gamma$ such that at
    least one of the ``noise driven'' vector fields is transverse
    to~$G'$.
  \end{enumerate}
\end{exmp}

%==============================================================================

\subsection{A remark concerning PDEs}

Notations can be deceiving, sometimes. The compact formulation
of~\eqref{eq:weak-FPK} and~\eqref{eq:FPK-spont}, which makes them look
very much like the usual Fokker-Planck equation, should not fool the
reader into thinking that these equations are simple PDEs. Indeed,
even when a (piecewise) smooth pdf exists, the generalized FPK
equation is in general a system of integro-differential equations,
with boundary conditions that can also involve integrals. The
integrals are hidden in the kernel notation: $(r_t K)(\Gamma) = \int
r_t(\dx) K(x,\Gamma)$. Fortunately, they disappear in many interesting
examples where the reset kernel is simple enough (see examples
\ref{ex:switching-diff}--\ref{ex:thermo}). This is an important
observation for practical applications, since the numerical solution
of a PDE is much easier than that of a general integro-differential
equation.

%==============================================================================

\appendix

\section{Proof of Theorem 4}
\label{sec:proof}

Let $\Ctest$ denote the set of all compactly supported~$\varphi \in
\Cscr^2(\Eset)$. The following lemma is an easy consequence of the
smoothness of the vector fields:

\pagebreak[2]

\begin{lem}
  For all $\varphi \in \Cscr^2(\Eset)$, $t \mapsto \int_0^t (\Gene^* \mu_s)(\varphi)\, \ds$ is
  differentiable on the right, with the right continuous derivative $t
  \mapsto (\Gene^* \mu_t)(\varphi)$.
\end{lem}

In the sequel, ``right continuous'' is abbreviated as ``rc''.

\newcommand{\dempoint}{$\;\diamond\;$}

\dempoint Assume that both~\Ha and~\Hb hold. Then each term
of~\eqref{eq:Dynkin-rewritten} has a $t$-derivative on the right.
Differentiating both sides proves that~\eqref{eq:weak-FPK} holds for
all~$t \geq 0$, hence that~$\Gene^* \mu_t$ is a Radon measure and that~$t
\mapsto \Gene^* \mu_t$ is rc. Moreover, integrating the inequality
$\left| \Gene^* \mu_t \right| \leq \left| \mu'_t \right| + 2 r_t$ yields
that, for all~$\Gamma \in \Ecalc$,
\begin{equation*}
  \int_0^t \left| \Gene^* \mu_s \right|(\Gamma)\, \ds
  \leq \int_0^t \left| \mu'_s \right| (\Gamma)\, \ds
  + 2\, \EE \bigl\{ N_t \bigr\} 
  \leq +\infty \,.
\end{equation*}
Therefore $t \mapsto \left| \Gene^* \mu_s \right|$ is locally integrable,
which proves~\Hc.

\dempoint Assume now that \Ha and \Hc hold, and set $\mu'_t = \Gene^*
\mu_t + r_t (K-I)$, for all~$t \geq 0$. Clearly, $\mu'_t$ is a Radon
measure, $t \mapsto \mu'_t$ is rc and
\begin{equation}
  \label{eq:demo:rel:int}
 \int_0^t \mu'_t \varphi = (\mu_t - \mu_0)\varphi \,,
 \quad \forall t \geq 0 \,,
 \quad \forall \varphi \in \Ctest \,.
\end{equation}
Moreover, for all~$\Gamma \in \Ecalc$,
\begin{equation*}
  \int_0^t \left| \mu'_s \right| (\Gamma)\, \ds
  \leq \int_0^t \left| \Gene^* \mu_s \right|(\Gamma)\, \ds +
  2\, \EE \bigl\{ N_t \bigr\} 
  \leq +\infty \,,
\end{equation*}
which shows that~$t \mapsto \left| \mu'_s \right|$ is locally
integrable. Therefore, using standard approximation techniques and a
monotone class argument, it can be proved that~\eqref{eq:demo:rel:int}
still holds for~$\varphi = \un_{\Gamma}$, $\Gamma \in \Ecalc$, i.e. that $t
\mapsto \mu'_t$ is the ``derivative'' of $t \mapsto \mu_t$ in the sense of
definition~\ref{def:GFPK:b}.

\dempoint Finally, assume that \Hb and \Hc hold.  Then, for
all~$\varphi \in \Ctest$, equation~$\eqref{eq:Dynkin-rewritten}$ can
be rewritten as
\begin{multline}
  \label{eq:zooiiing:2}
  \iint\limits_{G×]0;t]}\, \varphi(x)\, 
  \left(
    R^G(\dx,\ds) - (\Gene^* \mu_s)(\dx)\ds 
  \right) \hfill \\
  \;=\; 
  \iint\limits_{\EZset×]0;t]}\, \varphi(x)\, 
  \left(
    (R^G K)(\dx,\ds) - \xi_s(\dx)\,\ds 
  \right),
\end{multline}
where $\xi_s = \mu'_s - \bigl(\Gene^* \mu_s\bigr)(\EZset\cap\dotvar) -
r^0(K{-}I)$. The measures~$R^G$ and~$r^0$ have been defined in
subsection~\ref{subsec:def-and-link}. Clearly, $\xi_t \in \Mrad$ and $t
\mapsto \xi_t$ is locally integrable. Using once more standard
approximation techniques, one can prove that~\eqref{eq:zooiiing:2}
still holds when~$\varphi = \un_{\Gamma}$, with $\Gamma$ a compact subset
of~$G$. In this case the right-hand side vanishes, yielding
\begin{equation*}
  R^G(\Gamma × ]0;t]) 
  \;=\; \int_0^t (\Gene^* \mu_s)(\Gamma)\, \ds
  \,.
\end{equation*}
Moreover, since $t \mapsto R^G(\Gamma × ]0;t])$ is increasing and
$t \mapsto (\Gene^* \mu_t)(\Gamma)$ is rc, we have $(\Gene^* \mu_t)(\Gamma) \geq 0$
for all~$t \geq 0$. This allows to extend~\eqref{eq:zooiiing:2} to
all~$\Gamma \in \Ecalc$, using a monotone class argument, thus proving
the existence of a mean jump intensity $r_t^G = (\Gene^* \mu_s)(G \cap
\dotvar)$ for the forced jumps.

\hspace{1ex}

%==============================================================================

% Note: the \bibliographystyle command is issued in \AtBeginDocument
% by the document class 'infacconf.cls', it would be a mistake to do
% it again! (fucking template)

\bibliography{Bect-IFACWC2008-FinalSubmission}

\end{document}